\documentclass{article}

\usepackage{arxiv}

\usepackage[utf8]{inputenc} 
\usepackage[T1]{fontenc}    
\usepackage{hyperref}       
\usepackage{url}            
\usepackage{booktabs}       
\usepackage{amsfonts}       
\usepackage{nicefrac}       
\usepackage{microtype}      
\usepackage{lipsum}		
\usepackage{graphicx, amsmath}
\usepackage{doi}

\title{Nonlocal problem for multi-parametric integral-differential equation}

\author{ \href{https://orcid.org/0000-0003-4443-6300}{\includegraphics[scale=0.06]{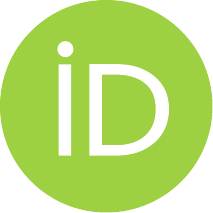}\hspace{1mm}Erkinjon Karimov}\thanks{EK is partially supported by the Methusalem programme of the Ghent University Special Research Fund (BOF), Grant/Award Number: 01M01021.} \\
	Department of Mathematics: Analysis, Logic and Discrete Mathematics\\
	Ghent University\\
	Ghent, Belgium \\
	\texttt{erkinjon.karimov@ugent.be} \\
	\And
	\href{https://orcid.org/0000-0002-3574-075X}{\includegraphics[scale=0.06]{orcid.pdf}\hspace{1mm}Doniyor Usmonov} \\
	Department of Mathematical Analysis and Differential Equations\\
	Fergana State University\\
	Fergana, Uzbekistan \\
	\texttt{dusmonov909@gmail.com} \\
	\And
	\href{https://orcid.org/0000-0003-4960-8669}{\includegraphics[scale=0.06]{orcid.pdf}\hspace{1mm}Khurshidjon Turdiev} \\
	Department of Mathematical Analysis and Differential Equations\\
	Fergana State University\\
	Fergana, Uzbekistan \\
	\texttt{xurshidjon2801@gmail.com} \\
}

\renewcommand{\shorttitle}{Nonlocal problem for multi-parametric integral-differential equation}

\hypersetup{
pdftitle={Nonlocal problem for multi-parametric integral-differential equation},
pdfauthor={Erkinjon Karimov, Doniyor Usmonov, Khurshidjon Turdiev},
pdfkeywords={Integral-differential equation, nonlocal problem, integral equation, Goursat problem, multi-parametric integral-differential operator, Mittag-Leffler type function.},
}

\begin{document}
\maketitle

\begin{abstract}
	This paper investigates a nonlocal boundary value problem for a multi-parametric integral-differential equation involving the Caputo-Prabhakar type operator in a bounded rectangular domain. The nonlocal conditions are given as partial integral expressions of the unknown function with continuous kernels. Using a known representation of the solution to the corresponding Goursat problem in terms of bivariate and trivariate Mittag-Leffler-type functions, the problem is reduced to a system of Volterra integral equations of the second kind for boundary traces. Based on this reduction, sufficient conditions ensuring existence and uniqueness of the solution are established. An explicit representation of the solution is also obtained via the solution of the derived integral system.  
\end{abstract}


\section*{Introduction}

Equations involving mixed partial derivatives of the form
\[
\frac{\partial^2 u(t,x)}{\partial t \partial x} - \lambda u(t,x) = f(t,x)
\]
play an important role in the theory of hyperbolic-type equations and arise naturally in the study of transport-type processes and related mathematical models. Such equations are closely connected with Goursat-type problems, where boundary conditions are prescribed along characteristic curves.

A powerful approach to the analysis of such problems is based on the use of integral representations of solutions. These representations allow one to reduce boundary value problems to systems of Volterra integral equations of the second kind, which provide an effective framework for establishing existence and uniqueness results.

A significant direction in the development of this theory concerns boundary conditions of nonlocal type. In such problems, the values of the unknown function on the boundary are connected through integral relations with given kernels. Problems with integral boundary conditions have been widely studied; for instance, in \cite{Psh20}, \cite{Psh_17}, \cite{Psh_16}, nonlocal conditions for related equations were considered, and the corresponding problems were reduced to Volterra integral equations. This approach made it possible to establish existence and uniqueness theorems and to obtain representations of solutions.

In recent years, significant attention has been devoted to the study of \textit{integral-differential equations with multi-parametric operators}. Among them, operators of Prabhakar-type (\cite{P71}) play an increasingly important role due to their flexibility and ability to incorporate memory effects through special kernels expressed in terms of Mittag-Leffler-type functions. These operators generalize classical fractional derivatives and allow the construction of more general analytical models.  Such operators allow one to construct multi-parametric integral-differential models that extend classical equations while preserving their analytical structure.

The analysis of these generalized equations is closely connected with the theory of Mittag-Leffler-type functions. The main results and methodological contributions of the works under consideration are discussed in detail in the subsequent paragraphs. These functions arise naturally in the representation of solutions to integral-differential equations and play a role analogous to the exponential function in classical analysis.

Using representations involving Mittag-Leffler-type functions, boundary value problems for multi-parametric integral-differential equations can be reduced to systems of Volterra integral equations. This reduction provides an effective tool for studying solvability, since the theory of integral equations guarantees existence and uniqueness of solutions under appropriate conditions on the kernels and given functions.

Despite the significant progress achieved in the study of nonlocal problems and integral-differential generalizations, equations involving multi-parametric operators of Prabhakar type combined with mixed derivative structures remain insufficiently investigated. In particular, the extension of equation \eqref{eq1} to the case of multi-parametric integral-differential operators leads to new classes of problems that require further analysis. This motivates the present study.

In parallel, substantial progress has been made in the theory of special functions tailored to fractional-order equations. Generalized Mittag-Leffler-type functions -- including multivariate and multiparameter extensions -- have emerged as natural solution kernels for fractional differential equations. A significant contribution in this area is the introduction of the bivariate Mittag-Leffler-type functions $E_{1}(x,y)$ and $E_{2}(x,y)$ by \cite{M. Garg}. Although the analytical properties of $E_{1}(x,y)$ were thoroughly investigated in that work, the corresponding properties of $E_{2}(x,y)$ remained unaddressed. This gap was subsequently filled by a series of studies \cite{E.T. Karimov}, \cite{E.T. Karimov1}, \cite{KKT26}, \cite{S. Kerbal}, which systematically examined the mathematical characteristics of $E_{2}(x,y)$ and demonstrated its applicability to fractional-order partial differential equations. In \cite{A.Hasanov}, a broader family of bivariate Mittag-Leffler-type functions, namely ${D_1}(x,y), \ldots, {D_5}(x,y)$ and ${E_3}(x,y), \ldots, {E_{11}}(x,y)$, was introduced, with particular emphasis placed on a detailed investigation of the properties of ${D_1}(x,y)$. In a similar vein, the function ${E_{12}}(x,y)$ was introduced and studied in \cite{TU25} and \cite{UM25}.

In \cite{HAGM20}, a novel trivariate Mittag-Leffler function was introduced, and its complex integral representation, Laplace transform relations, and fractional derivative properties were systematically investigated, along with its applications to the solution of multi-term fractional differential equations. The work in \cite{AHFM21} presents a new class of Mittag-Leffler functions defined through triple series, together with the development of associated fractional calculus operators, culminating in fully closed-form solutions. In \cite{HY24}, four new trivariate Mittag-Leffler-type functions $\overline{F}_{A}^{(3)}$, $\overline{F}_{B}^{(3)}$, $\overline{F}_{C}^{(3)}$, and $\overline{F}_{D}^{(3)}$ were introduced and thoroughly examined; these functions constitute natural generalizations of the classical Lauricella hypergeometric functions to three variables, and their fundamental properties, integral representations, and relevance to fractional differential equations were established in detail. Similarly, several further classes of trivariate Mittag-Leffler-type functions were introduced and investigated in \cite{KTU25}, \cite{KKh26}, \cite{RT25}.

Regarding boundary and nonlocal problems for fractional generalization of hyperbolic equations, a growing body of literature has emerged in recent years. R.A. Pshibikhova~\cite{Psh20} studied a nonlocal boundary value problem for a generalized hyperbolic equation with Caputo fractional derivatives in a bounded rectangular domain of the plane with two independent variables. The equation under consideration has the form
\begin{equation}\label{eq:Psh}
\partial_{0x}^{\alpha}\partial_{0y}^{\beta}u(x,y) + \lambda\, u(x,y) = f(x,y),
\end{equation}
subject to nonlocal boundary conditions prescribed as partial integral expressions of the unknown solution with respect to each variable:
\begin{equation}\label{eq:Psh_bc1}
    \int_{0}^{b} M(x,y)\,u(x,y)\,dy = \tau(x), \qquad 0 \leq x \leq a,
\end{equation}
\begin{equation}\label{eq:Psh_bc2}
    \int_{0}^{a} K(x,y)\,u(x,y)\,dx = \eta(y), \qquad 0 \leq y \leq b,
\end{equation}
where $M(x,y)$ and $K(x,y)$ are given continuous kernels. Employing a previously obtained representation of the solution to the Goursat problem for this equation in terms of Wright-type functions, the problem was reduced to a system of linear second-kind Volterra integral equations with respect to the traces of the unknown solution on a part of the domain boundary. As a result, the existence and uniqueness of the solution were established, and its representation was found in terms of the solutions to the obtained system of integral equations. In \cite{Psh_17}, \cite{Psh_16}, similar problems for equation \eqref{eq:Psh} were also investigated under mixed boundary-integral conditions. Kh.N. Turdiev \cite{T24} investigated two nonlocal boundary value problems for a system of coupled hyperbolic equations in a mixed pentagonal domain, deriving associated second-kind Volterra integral equations and proving solution uniqueness via the general theory of integral equations. Further contributions to this direction appear in \cite{AS26}, \cite{M17}.

Despite this progress, nonlocal problems for hyperbolic-type equations governed by the Caputo-Prabhakar operator -- a multi-parametric integral-differential operator generalizing both the classical Caputo derivative and the Prabhakar integral -- have received virtually no attention in the literature. This gap motivates the present work.

\textbf{Contribution of the paper.}\\
In this work, we investigate a nonlocal boundary value problem for a multi-parametric integral-differential equation in a bounded domain. The main contributions are:
\begin{itemize}
\item reduction of the problem to a system of Volterra integral equations;
\item proof of existence and uniqueness of the solution;
\item derivation of an explicit representation of the solution in terms of Mittag-Leffler-type functions.
\end{itemize}

\medskip

The obtained results contribute to the theory of nonlocal problems for integral-differential equations and extend existing methods based on integral equation techniques.

\section{Formulation of problem}

We consider the following multi-parametric integral-differential equation
\begin{equation}\label{eq1}
\frac{\partial }{\partial x}{{\,}^{PC}}D_{0t}^{\alpha ,\beta ,\gamma ,\delta }u\left( t,x \right)-a\frac{\partial }{\partial x}u\left( t,x \right)-b{{\,}^{PC}}D_{0t}^{\alpha ,\beta ,\gamma ,\delta }u\left( t,x \right)=f\left( t,x \right)
\end{equation}
in a domain $\Omega =\left\{ \left( t,x \right):0<t<q,0<x<p \right\}$. Here $f\left(t,x\right)$ is a given function and
\begin{equation} \label{GrindEQ__2_} 
\, ^{PC} D_{0t}^{\alpha ,\beta ,\gamma ,\delta } y\left(t\right)=\, \, ^{P} I_{0t}^{\alpha ,m-\beta ,-\gamma ,\delta } \frac{d^{m} }{dt^{m} } y\left(t\right) 
\end{equation} 
represents Caputo-Prabhakar multi-parametric integral-differential operator \cite{DP18} and 
\begin{equation} \label{GrindEQ__3_} 
\, \, ^{P} I_{0t}^{\alpha ,\beta ,\gamma ,\delta } y\left(t\right)=\int\limits _{0}^{t}\left(t-\xi \right)^{\beta -1} E_{\alpha ,\beta }^{\gamma } \left[\delta \left(t-\xi \right)^{\alpha } \right]y\left(\xi \right)d\xi ,\, \, \, \, \, t>0  
\end{equation} 
represents Prabhakar fractional integral \cite{P71}. We note that above-given definitions are valid for $\alpha,$ $\beta,$ $\gamma,$ $\delta,$ $a, b, q, p \in \mathbb{R}$ such that $\alpha >0$ and $m-1<\beta <m$.

Here $E_{\alpha ,\beta }^{\gamma } \left(z\right)$ is the generalized Mittag-Leffler function \cite{P71}:
$$E_{\alpha ,\beta }^{\gamma } \left(z\right)=\sum _{m=0}^{\infty }\frac{\left(\gamma \right)_{m} }{\Gamma \left(\alpha m+\beta \right)}  \frac{z^{m} }{m!}.$$

\textbf{Problem $N$.} We are interested to find a regular solution of the equation \eqref{eq1} with $0<\beta<1$ in $\Omega $, satisfying boundary condition
\begin{equation} \label{bc} 
u(t,0)=\varphi (t),\quad 0\le t\le q \end{equation} 
and nonlocal condition
\begin{equation} \label{nc} 
u(0,x)-\int\limits_{0}^{q}M(t)u(t,x)dt=\psi(x), \quad 0 \le x \le p,
\end{equation} 
where $\varphi \left( t \right)$, $\psi(x)$ and $M(t)$  are given functions ($M(t)\not\equiv 0$), such that $\varphi (0) -\int\limits_{0}^{q}M(t)\varphi(t)dt=\psi(0)$.

{\bf Definition.}
We call a function $u\left(t,x\right)$ as a regular solution of problem \eqref{eq1}, \eqref{bc}-\eqref{nc}, if
$u(t,x)\in C(\overline{\Omega}),\,\,u_{x}(t,x)\in C(\Omega),\,\, {{}^{PC}}D_{0t}^{\alpha ,\beta ,\gamma ,\delta }u(t,x)\in C(\Omega)$.

We note that at $\beta=1$, $\delta=0$, Eq. \eqref{eq1} becomes classical hyperbolic-type equation:
$$u_{tx}(t,x)-au_{x}(t,x)-bu_{t}(t,x)=f(t,x).$$

\section{Main results}

Let us introduce a notation  
\begin{equation}\label{ic}
u(0,x)=\tau(x), \quad 0 \le x \le p.    
\end{equation}

{\bf Theorem 1. }
If $\varphi (0)=\tau (0)$, $\tau(x)\in C[0,1]\cap C^{1}(0,1)$ and $\varphi(t)\in C[0,q]\cap C^1(0,q)$, then the solution of Goursat problem multi-parametric integral-differential equation \eqref{eq1}, that satisfies the initial condition \eqref{ic} and boundary condition \eqref{bc} will be represented as follows:
$$u(t,x)=\tau \left( x \right)+\left( \varphi \left( t \right)-\varphi \left( 0 \right) \right){{e}^{bx}}+a\Gamma \left( \gamma  \right){{t}^{\beta }}\tau \left( x \right){{E}_{2}}\left( \left. \begin{aligned}
  & \gamma ,1,\gamma ;0,1; \\ 
 & \beta ,\alpha ,\beta +1;\gamma ,\gamma ;1,1; \\ 
\end{aligned} \right|\begin{matrix}
   a{{t}^{\beta }}  \\
   \delta {{t}^{\alpha }}  \\
\end{matrix} \right)-$$
$$a\varphi \left( 0 \right){{t}^{\beta }}\overline{F}_{E}^{(3)}\left( \left. \begin{aligned}
  & \gamma ,1,\gamma ;1,1,2; \\ 
 & \beta ,\alpha ,\beta +1;\gamma ,\gamma ;1,2;1,1;1,1;1,1; \\ 
\end{aligned} \right|a{{t}^{\beta }};bx;\delta {{t}^{\alpha }} \right)+$$
$$ab{{t}^{\beta }}\int\limits_{0}^{x}{\tau \left( \xi  \right)\overline{F}_{E}^{(3)}\left( \left. \begin{aligned}
  & \gamma ,1,\gamma ;1,1,2; \\ 
 & \beta ,\alpha ,\beta +1;\gamma ,\gamma ;1,1;1,1;1,1;1,2; \\ 
\end{aligned} \right|a{{t}^{\beta }};b\left( x-\xi  \right);\delta {{t}^{\alpha }} \right)d\xi }+$$
$$abx\int\limits_{0}^{t}{{{\left( t-\eta  \right)}^{\beta -1}}\varphi \left( \eta  \right)\overline{F}_{E}^{(3)}\left( \left. \begin{aligned}
  & \gamma ,1,\gamma ;1,1,2; \\ 
 & \beta ,\alpha ,\beta ;\gamma ,\gamma ;1,2;1,1;1,1;1,2; \\ 
\end{aligned} \right|a{{\left( t-\eta  \right)}^{\beta }};bx;\delta {{\left( t-\eta  \right)}^{\alpha }} \right)d\eta }+$$
\begin{equation}\label{u(t,x)}
  \int\limits_{0}^{t}{\int\limits_{0}^{x}{{{\left( t-\eta  \right)}^{\beta -1}}f\left( \eta ,\xi  \right)\overline{F}_{E}^{(3)}\left( \left. \begin{aligned}
  & \gamma ,1,\gamma ;1,1,1; \\ 
 & \beta ,\alpha ,\beta ;\gamma ,\gamma ;1,1;1,1;1,1;1,1; \\ 
\end{aligned} \right|a{{\left( t-\eta  \right)}^{\beta }};b\left( x-\xi  \right);\delta {{\left( t-\eta  \right)}^{\alpha }} \right)d\eta d\xi }}.
\end{equation}
The solution to the Goursat problem exists and is unique (see for details \cite{KTU25}).

Here, $E_{2}(\cdot)$ the bivariate Mittag-Leffler type function \cite{M. Garg}:
$${{E}_{2}}\left( \left. \begin{aligned}
  & {{\alpha }_{1}},{{\beta }_{1}},{{\gamma }_{1}};{{\alpha }_{2}},{{\gamma }_{2}}; \\ 
 & {{\alpha }_{3}},{{\beta }_{2}},{{\delta }_{1}};{{\alpha }_{4}},{{\delta }_{2}};{{\beta }_{3}},{{\delta }_{3}}; \\ 
\end{aligned} \right|\begin{matrix}
   x  \\
   y  \\
\end{matrix} \right)=$$
\begin{equation}\label{E2qat}
\sum\limits_{m=0}^{\infty }{\sum\limits_{k=0}^{\infty }{\frac{\Gamma \left( {{\alpha }_{1}}m+{{\beta }_{1}}k+{{\gamma }_{1}} \right)\Gamma \left( {{\alpha }_{2}}m+{{\gamma }_{2}} \right)}{\Gamma \left( {{\gamma }_{1}} \right)\Gamma \left( {{\gamma }_{2}} \right)\Gamma \left( {{\alpha }_{3}}m+{{\beta }_{2}}k+{{\delta }_{1}} \right)}\frac{{{x}^{m}}}{\Gamma \left( {{\alpha }_{4}}m+{{\delta }_{2}} \right)}\frac{{{y}^{k}}}{\Gamma \left( {{\beta }_{3}}k+{{\delta }_{3}} \right)}}}, 
\end{equation} 
$$\gamma_{1}, \gamma_{2}, \delta_{1}, \delta_{2}, \delta_{3}, x,y\in \mathbb{R}, \,\, \min \left\{ \alpha_{1}, \alpha_{2}, \alpha_{3}, \alpha_{4}, \beta_{1}, \beta_{2}, \beta_{3} \right\}>0,$$
in which the double series converges for $x,y\in \mathbb{R}$, if ${{\alpha }_{3}}+{{\alpha }_{4}}-{{\alpha }_{1}}-{{\alpha }_{2}}>0$, and ${{\beta }_{2}}+{{\beta }_{3}}-{{\beta }_{1}}>0$. 

$\overline{F}_{E}^{(3)}(\cdot)$ the trivariate Mittag-Leffler type function \cite{KTU25}:
$$\overline{F}_{E}^{(3)}\left( \left. \begin{aligned}
  & {{\alpha }_{1}},{{\beta }_{1}},{{\delta }_{1}};{{\alpha }_{2}},{{\gamma }_{1}},{{\delta }_{2}}; \\ 
 & {{\alpha }_{3}},{{\beta }_{2}},{{\delta }_{3}};{{\alpha }_{4}},{{\delta }_{4}};{{\alpha }_{5}},{{\delta }_{5}};{{\beta }_{3}},{{\delta }_{6}};{{\gamma }_{2}},{{\delta }_{7}};{{\gamma }_{3}},{{\delta }_{8}}; \\ 
\end{aligned} \right|x;y;z \right)=$$
\begin{equation}\label{Hs}
       \sum\limits_{m,j,k=0}^{+\infty }{\frac{\Gamma \left( {{\alpha }_{1}}m+{{\beta }_{1}}k+{{\delta }_{1}} \right)\Gamma \left( {{\alpha }_{2}}m+{{\gamma }_{1}}j+{{\delta }_{2}} \right){{x}^{m}}{{y}^{j}}{{z}^{k}}}{\Gamma \left( {{\alpha }_{3}}m+{{\beta }_{2}}k+{{\delta }_{3}} \right)\Gamma \left( {{\alpha }_{4}}m+{{\delta }_{4}} \right)\Gamma \left( {{\alpha }_{5}}m+{{\delta }_{5}} \right)\Gamma \left( {{\beta }_{3}}k+{{\delta }_{6}} \right)\Gamma \left( {{\gamma }_{2}}j+{{\delta }_{7}} \right)\Gamma \left( {{\gamma }_{3}}j+{{\delta }_{8}} \right)}},
\end{equation}
$$\left({{\alpha }_{l}},{{\beta }_{i}},{{\gamma }_{i}}, {{\delta }_{n}},x,y,z\in \mathbb{R};\,\,\,\min \left\{ {{\alpha }_{l}},{{\beta }_{i}},{{\gamma }_{i}} \right\}>0,\,\,\,\left( n=\overline{1,8},\,\,\,l=\overline{1,5},\,\,\,i=\overline{1,3} \right) \right),$$
in which the triple series converges for $x,y,z\in \mathbb{R}$, if $\Delta_{1}>0$, $\Delta_{2}>0$ and $\Delta_{3}>0$. Here
$${{\Delta }_{1}}={{\alpha }_{3}}+{{\alpha }_{4}}+{{\alpha }_{5}}-{{\alpha }_{1}}-{{\alpha }_{2}}, \,\,\, {{\Delta }_{2}}={{\gamma }_{2}}+{{\gamma }_{3}}-{{\gamma }_{1}},  \,\,\, {{\Delta }_{3}}={{\beta }_{2}}+{{\beta }_{3}}-{{\beta }_{1}}.$$

To determine the function $\tau(x)$, we will use condition \eqref{nc}. Following the representation \eqref{u(t,x)}, and considering the nonlocal condition \eqref{nc}, we obtain 
$$\tau \left( x \right)-\tau \left( x \right)\int\limits_{0}^{q}{M}\left( t \right)dt-{{e}^{bx}}\int\limits_{0}^{q}{M}\left( t \right)\left( \varphi \left( t \right)-\varphi \left( 0 \right) \right)dt-$$
$$a\Gamma(\gamma)\tau(x)\int\limits_{0}^{q}M(t){{t}^{\beta }}{{E}_{2}}\left( \left. \begin{aligned}
  & \gamma ,1,\gamma ;0,1; \\ 
 & \beta ,\alpha ,\beta +1;\gamma ,\gamma ;1,1; \\ 
\end{aligned} \right|\begin{matrix}
   a{{t}^{\beta }}  \\
   \delta {{t}^{\alpha }}  \\
\end{matrix} \right)dt+$$
$$a\varphi(0)\int\limits_{0}^{q}M(t){{t}^{\beta }}\overline{F}_{E}^{(3)}\left( \left. \begin{aligned}
  & \gamma ,1,\gamma ;1,1,2; \\ 
 & \beta ,\alpha ,\beta +1;\gamma ,\gamma ;1,2;1,1;1,1;1,1; \\ 
\end{aligned} \right|a{{t}^{\beta }};bx;\delta {{t}^{\alpha }} \right)dt-$$
$$ab\int\limits_{0}^{q}M(t)\int\limits_{0}^{x}{\tau \left( \xi  \right){{t}^{\beta }}\overline{F}_{E}^{(3)}\left( \left. \begin{aligned}
  & \gamma ,1,\gamma ;1,1,2; \\ 
 & \beta ,\alpha ,\beta +1;\gamma ,\gamma ;1,1;1,1;1,1;1,2; \\ 
\end{aligned} \right|a{{t}^{\beta }};b\left( x-\xi  \right);\delta {{t}^{\alpha }} \right)d\xi }dt-$$
$$abx\int\limits_{0}^{q}M(t)\int\limits_{0}^{t}{{\left( t-\eta  \right)}^{\beta -1}}\varphi(\eta)\times$$
$$\overline{F}_{E}^{(3)}\left( \left. \begin{aligned}
  & \gamma ,1,\gamma ;1,1,2; \\ 
 & \beta ,\alpha ,\beta ;\gamma ,\gamma ;1,2;1,1;1,1;1,2; \\ 
\end{aligned} \right|a{{\left( t-\eta  \right)}^{\beta }};bx;\delta {{\left( t-\eta  \right)}^{\alpha }} \right)dtd\eta-$$
$$\int\limits_{0}^{q}M(t)\int\limits_{0}^{t}\int\limits_{0}^{x}{{\left( t-\eta  \right)}^{\beta -1}}f\left( \eta ,\xi  \right)\times$$
\begin{equation}\label{u(q,0)}
\overline{F}_{E}^{(3)}\left( \left. \begin{aligned}
  & \gamma ,1,\gamma ;1,1,1; \\ 
 & \beta ,\alpha ,\beta ;\gamma ,\gamma ;1,1;1,1;1,1;1,1; \\ 
\end{aligned} \right|a{{\left( t-\eta  \right)}^{\beta }};b\left( x-\xi  \right);\delta {{\left( t-\eta  \right)}^{\alpha }} \right)d\eta d\xi dt=\psi(x).
\end{equation}

By changing the order of integration, from equation \eqref{u(q,0)}, we deduce
$$\tau \left( x \right)-\tau \left( x \right)\int\limits_{0}^{q}{M}\left( t \right)dt-{{e}^{bx}}\int\limits_{0}^{q}{M}\left( t \right)\left( \varphi \left( t \right)-\varphi \left( 0 \right) \right)dt-$$
$$a\Gamma(\gamma)\tau(x)\int\limits_{0}^{q}M(t){{t}^{\beta }}{{E}_{2}}\left( \left. \begin{aligned}
  & \gamma ,1,\gamma ;0,1; \\ 
 & \beta ,\alpha ,\beta +1;\gamma ,\gamma ;1,1; \\ 
\end{aligned} \right|\begin{matrix}
   a{{t}^{\beta }}  \\
   \delta {{t}^{\alpha }}  \\
\end{matrix} \right)dt+$$
$$a\varphi(0)\int\limits_{0}^{q}M(t){{t}^{\beta }}\overline{F}_{E}^{(3)}\left( \left. \begin{aligned}
  & \gamma ,1,\gamma ;1,1,2; \\ 
 & \beta ,\alpha ,\beta +1;\gamma ,\gamma ;1,2;1,1;1,1;1,1; \\ 
\end{aligned} \right|a{{t}^{\beta }};bx;\delta {{t}^{\alpha }} \right)dt-$$
$$ab\int\limits_{0}^{x}\tau \left( \xi  \right)\int\limits_{0}^{q}M(t){{{t}^{\beta }}\overline{F}_{E}^{(3)}\left( \left. \begin{aligned}
  & \gamma ,1,\gamma ;1,1,2; \\ 
 & \beta ,\alpha ,\beta +1;\gamma ,\gamma ;1,1;1,1;1,1;1,2; \\ 
\end{aligned} \right|a{{t}^{\beta }};b\left( x-\xi  \right);\delta {{t}^{\alpha }} \right)dtd\xi }-$$
$$abx\int\limits_{0}^{q}\varphi(\eta)\int\limits_{\eta}^{q}M(t){{\left( t-\eta  \right)}^{\beta -1}}\times$$
$$\overline{F}_{E}^{(3)}\left( \left. \begin{aligned}
  & \gamma ,1,\gamma ;1,1,2; \\ 
 & \beta ,\alpha ,\beta ;\gamma ,\gamma ;1,2;1,1;1,1;1,2; \\ 
\end{aligned} \right|a{{\left( t-\eta  \right)}^{\beta }};bx;\delta {{\left( t-\eta  \right)}^{\alpha }} \right)dtd\eta-$$
$$\int\limits_{0}^{q}M(t)\int\limits_{0}^{t}\int\limits_{0}^{x}{{\left( t-\eta  \right)}^{\beta -1}}f\left( \eta ,\xi  \right)\times$$
\begin{equation}\label{u(q,0)1}
\overline{F}_{E}^{(3)}\left( \left. \begin{aligned}
  & \gamma ,1,\gamma ;1,1,1; \\ 
 & \beta ,\alpha ,\beta ;\gamma ,\gamma ;1,1;1,1;1,1;1,1; \\ 
\end{aligned} \right|a{{\left( t-\eta  \right)}^{\beta }};b\left( x-\xi  \right);\delta {{\left( t-\eta  \right)}^{\alpha }} \right)d\eta d\xi dt=\psi(x).
\end{equation}

Let us introduce the following notation:
$$A=1-\int\limits_{0}^{q}{M}\left( t \right)\left( 1+a\Gamma \left( \gamma  \right){{t}^{\beta }}{{E}_{2}}\left( \left. \begin{aligned}
  & \gamma ,1,\gamma ;0,1; \\ 
 & \beta ,\alpha ,\beta +1;\gamma ,\gamma ;1,1; \\ 
\end{aligned} \right|\begin{matrix}
   a{{t}^{\beta }}  \\
   \delta {{t}^{\alpha }}  \\
\end{matrix} \right) \right)dt,$$
$$M_{1}(\xi,x)=\int\limits_{0}^{q}M(t){{{t}^{\beta }}\overline{F}_{E}^{(3)}\left( \left. \begin{aligned}
  & \gamma ,1,\gamma ;1,1,2; \\ 
 & \beta ,\alpha ,\beta +1;\gamma ,\gamma ;1,1;1,1;1,1;1,2; \\ 
\end{aligned} \right|a{{t}^{\beta }};b\left( x-\xi  \right);\delta {{t}^{\alpha }} \right)dt},$$
$$g(x)= \psi(x)+{{e}^{bx}}\int\limits_{0}^{q}M(t)\left( \varphi(t)-\varphi(0) \right)dt-$$
$$a\varphi(0)\int\limits_{0}^{q}M(t){{t}^{\beta }}\overline{F}_{E}^{(3)}\left( \left. \begin{aligned}
  & \gamma ,1,\gamma ;1,1,2; \\ 
 & \beta ,\alpha ,\beta +1;\gamma ,\gamma ;1,2;1,1;1,1;1,1; \\ 
\end{aligned} \right|a{{t}^{\beta }};bx;\delta {{t}^{\alpha }} \right)dt+$$
$$abx\int\limits_{0}^{q}\varphi(\eta)\int\limits_{\eta}^{q}M(t){{\left( t-\eta  \right)}^{\beta -1}}\times$$
$$\overline{F}_{E}^{(3)}\left( \left. \begin{aligned}
  & \gamma ,1,\gamma ;1,1,2; \\ 
 & \beta ,\alpha ,\beta ;\gamma ,\gamma ;1,2;1,1;1,1;1,2; \\ 
\end{aligned} \right|a{{\left( t-\eta  \right)}^{\beta }};bx;\delta {{\left( t-\eta  \right)}^{\alpha }} \right)dtd\eta+$$
$$\int\limits_{0}^{q}M(t)\int\limits_{0}^{t}\int\limits_{0}^{x}{{\left( t-\eta  \right)}^{\beta -1}}f\left( \eta ,\xi  \right)\times$$
$$\overline{F}_{E}^{(3)}\left( \left. \begin{aligned}
  & \gamma ,1,\gamma ;1,1,1; \\ 
 & \beta ,\alpha ,\beta ;\gamma ,\gamma ;1,1;1,1;1,1;1,1; \\ 
\end{aligned} \right|a{{\left( t-\eta  \right)}^{\beta }};b\left( x-\xi  \right);\delta {{\left( t-\eta  \right)}^{\alpha }} \right)d\eta d\xi dt.$$

Considering the introduced notations, from equation \eqref{u(q,0)1} we deduce 
\begin{equation}\label{(3.2.4)}
A\tau(x)-ab\int\limits_{0}^{x}\tau(\xi)M_{1}(\xi,x)d\xi=g(x).    
\end{equation}

\textbf{Case 1.} Let $A\ne 0$. Then, denoting
$$M_{2}(\xi,x)=\frac{M_{1}(\xi,x)}{A}, \quad G(x)=\frac{g(x)}{A},$$ 
from the last relation, we deduce
\begin{equation}\label{ie}
  \tau(x)-ab\int\limits_{0}^{x}\tau(\xi)M_{2}(\xi,x)d\xi=G(x).  
\end{equation}

{\bf Theorem 2. }
Let $A\ne 0$, $a<0$, $b<0$, $\delta<0$. If $\varphi(0)=G(0)$, $\varphi(t)\in C[0,q]\cap C^1(0,q)$, $\psi(x)\in C[0,q]\cap C^1(0,q)$, $M(t)\in C[0,q]$, $f\left( t,x \right)={{t}^{-{{\varepsilon }_{1}}}}{{x}^{-{{\varepsilon }_{2}}}}{{\widetilde{f}}_{1}}\left( t,x \right)$, ${{\widetilde{f}}_{1}}\left( t,x \right)\in C\left( \overline{\Omega } \right)$ and $0\le {{\varepsilon }_{1}}<\beta$, $0\le {{\varepsilon }_{2}}<1$, then the solution to the problem $N$ exists and is unique.

Here, $G(0)=g(0)/A$,
$$g(0)=\varphi(0)- \varphi(0)\int\limits_{0}^{q}M(t)dt-a\Gamma(\gamma)\varphi(0)\int\limits_{0}^{q}M(t){{t}^{\beta }}{{E}_{2}}\left( \left. \begin{aligned}
  & \gamma ,1,\gamma ;0,1; \\ 
 & \beta ,\alpha ,\beta +1;\gamma ,\gamma ;1,1; \\ 
\end{aligned} \right|\begin{matrix}
   a{{t}^{\beta }}  \\
   \delta {{t}^{\alpha }}  \\
\end{matrix} \right)dt.$$

{\bf Proof. }
According to the theory of integral equations, Volterra's second-kind integral equation exists and is unique, if the kernel ${{M}_{2}}\left( \xi ,x \right)$ of the integral equation and the right-hand side function $G\left( x \right)$ are continuous.

For the integral equation (12) to have a continuous $M_{2}(\xi,x)$ kernel, the $M_{1}(\xi,x)$ kernel must be continuous. Let $M(t)\in C[0,q]$. Then, using the boundedness property of the function $\overline{F}_{E}^{(3)}$ (see Lemma 3.3, \cite{KTU25}) for negative arguments, we obtain:
$$M_{1}(\xi,x)=\int\limits_{0}^{q}M(t){{t}^{\beta }}\overline{F}_{E}^{(3)}\left( \left. \begin{aligned}
  & \gamma ,1,\gamma ;1,1,2; \\ 
 & \beta ,\alpha ,\beta +1;\gamma ,\gamma ;1,1;1,1;1,1;1,2; \\ 
\end{aligned} \right|a{{t}^{\beta }};b\left( x-\xi  \right);\delta {{t}^{\alpha }} \right)dt\le$$
$$C_{1}\int\limits_{0}^{q}M(t){t}^{\beta }dt\le C_{2},$$
where $C_1,C_2>0$ are constants. Hence, the kernel $M_{1}\left(\xi,x\right)$ is continuous.

If $g(x)$ is continuous, then $G(x)$ will be also continuous. If $\psi(x)\in C[0,p]$, $\varphi(t)\in C[0,q]$, $M(t)\in C[0,q]$ and $f\left( t,x \right)={{t}^{-{{\varepsilon }_{1}}}}{{x}^{-{{\varepsilon }_{2}}}}{{\widetilde{f}}_{1}}\left( t,x \right)$, ${{\widetilde{f}}_{1}}\left( t,x \right)\in C\left( \overline{\Omega } \right)$ and $0\le {{\varepsilon }_{1}}<\beta$, $0\le {{\varepsilon }_{2}}<1$, then we could prove the continuity of the function $g(x)$ using the boundedness property of the function $\overline{F}_{E}^{(3)}$ (see Lemma 3.3, \cite{KTU25}):
$$g(x)= \psi(x)+{{e}^{bx}}\int\limits_{0}^{q}M(t)\left( \varphi(t)-\varphi(0) \right)dt-$$
$$a\varphi(0)\int\limits_{0}^{q}M(t){{t}^{\beta }}\overline{F}_{E}^{(3)}\left( \left. \begin{aligned}
  & \gamma ,1,\gamma ;1,1,2; \\ 
 & \beta ,\alpha ,\beta +1;\gamma ,\gamma ;1,2;1,1;1,1;1,1; \\ 
\end{aligned} \right|a{{t}^{\beta }};bx;\delta {{t}^{\alpha }} \right)dt+$$
$$abx\int\limits_{0}^{q}\varphi(\eta)\int\limits_{\eta}^{q}M(t){{\left( t-\eta  \right)}^{\beta -1}}\times$$
$$\overline{F}_{E}^{(3)}\left( \left. \begin{aligned}
  & \gamma ,1,\gamma ;1,1,2; \\ 
 & \beta ,\alpha ,\beta ;\gamma ,\gamma ;1,2;1,1;1,1;1,2; \\ 
\end{aligned} \right|a{{\left( t-\eta  \right)}^{\beta }};bx;\delta {{\left( t-\eta  \right)}^{\alpha }} \right)dtd\eta+$$
$$\int\limits_{0}^{q}M(t)\int\limits_{0}^{t}\int\limits_{0}^{x}{{\left( t-\eta  \right)}^{\beta -1}}f\left( \eta ,\xi  \right)\times$$
$$\overline{F}_{E}^{(3)}\left( \left. \begin{aligned}
  & \gamma ,1,\gamma ;1,1,1; \\ 
 & \beta ,\alpha ,\beta ;\gamma ,\gamma ;1,1;1,1;1,1;1,1; \\ 
\end{aligned} \right|a{{\left( t-\eta  \right)}^{\beta }};b\left( x-\xi  \right);\delta {{\left( t-\eta  \right)}^{\alpha }} \right)d\eta d\xi dt\le C_{3}.$$

\textbf{Case 2.} 
If $A=0$, then the integral equation (3.2.4) will become the first-kind of Volterra equation 
\begin{equation}\label{(3.2.6)}
-ab\int\limits_{0}^{x}{\tau }\left( \xi  \right){{M}_{1}}\left( \xi ,x \right)d\xi =g\left( x \right).  
\end{equation}
By taking the once derivative of the integral equation \eqref{(3.2.6)}, we then obtain the following integral equation:
$$\tau \left( x \right)+\int\limits_{0}^{x}{\tau }\left( \xi  \right){{M}^{*}}\left( \xi ,x \right)d\xi ={{g}^{*}}\left( x \right).$$
Here
$${{M}^{*}}\left( \xi ,x \right)=\frac{{{\left[ {{M}_{1}}\left( \xi ,x \right) \right]}_{x}}}{B}, \quad {{g}^{*}}\left( x \right)=-\frac{1}{ab}\frac{{g}'\left( x \right)}{B},$$
$${{\left[ {{M}_{1}}\left( \xi ,x \right) \right]}_{x}}=\int\limits_{0}^{q}{M}\left( t \right){{t}^{\beta }} \overline{F}_{E}^{(3)}\left( \left. \begin{aligned}
  & \gamma ,1,\gamma ;1,1,3; \\ 
 & \beta ,\alpha ,\beta +1;\gamma ,\gamma ;1,1;1,1;1,1;1,3; \\ 
\end{aligned} \right|a{{t}^{\beta }};b\left( x-\xi  \right);\delta {{t}^{\alpha }} \right)dt,$$
$${g}'\left( x \right)={\psi }'\left( x \right)+b{{e}^{bx}}\int\limits_{0}^{q}{M}\left( t \right)\left( \varphi \left( t \right)-\varphi \left( 0 \right) \right)dt-$$
$$ab\varphi \left( 0 \right)\int\limits_{0}^{q}{M}\left( t \right){{t}^{\beta }}\overline{F}_{E}^{\left( 3 \right)}\left( \left. \begin{aligned}
  & \gamma ,1,\gamma ;1,1,3; \\ 
 & \beta ,\alpha ,\beta +1;\gamma ,\gamma ;1,2;1,1;1,1;1,2; \\ 
\end{aligned} \right|a{{t}^{\beta }};bx;\delta {{t}^{\alpha }} \right)dt+$$
$$ab\int\limits_{0}^{q}{\varphi }\left( \eta  \right)\int\limits_{\eta }^{q}{M}\left( t \right){{\left( t-\eta  \right)}^{\beta -1}}\times $$
$$\overline{F}_{E}^{\left( 3 \right)}\left( \left. \begin{aligned}
  & \gamma ,1,\gamma ;1,1,2; \\ 
 & \beta ,\alpha ,\beta ;\gamma ,\gamma ;1,2;1,1;1,1;1,1; \\ 
\end{aligned} \right|a{{\left( t-\eta  \right)}^{\beta }};bx;\delta {{\left( t-\eta  \right)}^{\alpha }} \right)dtd\eta +$$
$$\Gamma \left( \gamma  \right)\int\limits_{0}^{q}{f\left( \eta ,x \right)}\int\limits_{\eta }^{q}{M\left( t \right){{\left( t-\eta  \right)}^{\beta -1}}{{E}_{2}}\left( \left. \begin{aligned}
  & \gamma ,1,\gamma ;0,1; \\ 
 & \beta ,\alpha ,\beta ;\gamma ,\gamma ;1,1; \\ 
\end{aligned} \right|\begin{matrix}
   a{{\left( t-\eta  \right)}^{\beta }}  \\
   \delta {{\left( t-\eta  \right)}^{\alpha }}  \\
\end{matrix} \right)d\eta dt+}$$
$$b\int\limits_{0}^{q}{M}\left( t \right)\int\limits_{0}^{t}{\int\limits_{0}^{x}{{{\left( t-\eta  \right)}^{\beta -1}}}}f\left( \eta ,\xi  \right)\times $$
$$\overline{F}_{E}^{\left( 3 \right)}\left( \left. \begin{aligned}
  & \gamma ,1,\gamma ;1,1,2; \\ 
 & \beta ,\alpha ,\beta ;\gamma ,\gamma ;1,1;1,1;1,1;1,2; \\ 
\end{aligned} \right|a{{\left( t-\eta  \right)}^{\beta }};b\left( x-\xi  \right);\delta {{\left( t-\eta  \right)}^{\alpha }} \right)d\eta d\xi dt,$$
$$B={{M}_{1}}\left( x,x \right)=\Gamma \left( \gamma  \right)\int\limits_{0}^{q}{M}\left( t \right){{t}^{\beta }}{{E}_{2}}\left( \left. \begin{aligned}
  & \gamma ,1,\gamma ;0,1; \\ 
 & \beta ,\alpha ,\beta ;\gamma ,\gamma ;1,1; \\ 
\end{aligned} \right|\begin{matrix}
   a{{t}^{\beta }}  \\
   \delta {{t}^{\alpha }}  \\
\end{matrix} \right)dt.$$

{\bf Theorem 3. }
Let $A=0$, $\alpha =1$, $0<\beta <1$, $\gamma =\beta$, $a<0$, $b<0$, $\delta <0$. If $\tau \left( 0 \right)={{g}^{*}}\left( 0 \right)$, $\varphi \left( t \right)\in C\left[ 0,q \right]\cap {{C}^{1}}\left( 0,q \right)$, $\psi \left( x \right)\in {{C}^{1}}\left[ 0,q \right]\cap {{C}^{2}}\left( 0,q \right)$, $M\left( t \right)\in C\left[ 0,q \right]$, $f\left( t,x \right)={{t}^{-{{\varepsilon }_{1}}}}{{x}^{-{{\varepsilon }_{2}}}}\widetilde{{{f}_{1}}}\left( t,x \right)$, $\widetilde{{{f}_{1}}}\left( t,x \right)\in C\left( {\bar{\Omega }} \right)$ and $0\le {{\varepsilon }_{1}}<\beta$, $0\le {{\varepsilon }_{2}}<1$, then the solution to Problem $N$ exists and is unique.

Here ${{g}^{*}}\left( 0 \right)=-\left[ {g}'\left( 0 \right) \right]/\left[ abB \right]$,
$${g}'\left( 0 \right)={\psi }'\left( 0 \right)+b\int\limits_{0}^{q}{M}\left( t \right)\left( \varphi \left( t \right)-\varphi \left( 0 \right) \right)dt-$$
$$ab\varphi \left( 0 \right)\int\limits_{0}^{q}{M}\left( t \right){{t}^{\beta }}{{E}_{13}}\left( \left. \begin{aligned}
  & \gamma ,1,\gamma ;1,3; \\ 
 & \beta ,\alpha ,\beta +1;\gamma ,\gamma ;1,2;1,1; \\ 
\end{aligned} \right|\begin{matrix}
   a{{t}^{\beta }}  \\
   \delta {{t}^{\alpha }}  \\
\end{matrix} \right)dt+$$
$$ab\int\limits_{0}^{q}{\varphi }\left( \eta  \right)\int\limits_{\eta }^{q}{M}\left( t \right){{\left( t-\eta  \right)}^{\beta -1}}{{E}_{2}}\left( \left. \begin{aligned}
  & \gamma ,1,\gamma ;0,1; \\ 
 & \beta ,\alpha ,\beta ;\gamma ,\gamma ;1,1; \\ 
\end{aligned} \right|\begin{matrix}
   a{{\left( t-\eta  \right)}^{\beta }}  \\
   \delta {{\left( t-\eta  \right)}^{\alpha }}  \\
\end{matrix} \right)dtd\eta +$$
$$\Gamma \left( \gamma  \right)\int\limits_{0}^{q}{f\left( \eta ,x \right)}\int\limits_{\eta }^{q}{M\left( t \right){{\left( t-\eta  \right)}^{\beta -1}}{{E}_{2}}\left( \left. \begin{aligned}
  & \gamma ,1,\gamma ;0,1; \\ 
 & \beta ,\alpha ,\beta ;\gamma ,\gamma ;1,1; \\ 
\end{aligned} \right|\begin{matrix}
   a{{\left( t-\eta  \right)}^{\beta }}  \\
   \delta {{\left( t-\eta  \right)}^{\alpha }}  \\
\end{matrix} \right)d\eta dt},$$
where
$${{E}_{13}}\left( \left. \begin{aligned}
  & {{\alpha }_{1}},{{\beta }_{1}},{{\delta }_{1}};{{\alpha }_{2}},{{\delta }_{2}}; \\ 
 & {{\alpha }_{3}},{{\beta }_{2}},{{\delta }_{3}};{{\alpha }_{4}},{{\delta }_{4}};{{\alpha }_{5}},{{\delta }_{5}};{{\beta }_{3}},{{\delta }_{6}}; \\ 
\end{aligned} \right|\begin{matrix}
   x  \\
   y  \\
\end{matrix} \right) =$$  
$$\sum\limits_{m,k=0}^{+\infty }{\frac{\Gamma \left( {{\alpha }_{1}}m+{{\beta }_{1}}k+{{\delta }_{1}} \right)\Gamma \left( {{\alpha }_{2}}m+{{\delta }_{2}} \right){{x}^{m}}{{y}^{k}}}{\Gamma \left( {{\alpha }_{3}}m+{{\beta }_{2}}k+{{\delta }_{3}} \right)\Gamma \left( {{\alpha }_{4}}m+{{\delta }_{4}} \right)\Gamma \left( {{\alpha }_{5}}m+{{\delta }_{5}} \right)\Gamma \left( {{\beta }_{3}}k+{{\delta }_{6}} \right)}}, $$
$$\left({{\alpha }_{l}},{{\beta }_{i}}, {{\delta }_{n}},x,y\in \mathbb{R};\, \min \left\{ {{\alpha }_{l}},{{\beta }_{i}} \right\}>0;\, \left( n=\overline{1,6}, \, l=\overline{1,5}, \, i=\overline{1,3} \right) \right),$$
in which the double series converges for $x,y\in \mathbb{R}$, if ${{\alpha }_{3}}+{{\alpha }_{4}}+{{\alpha }_{5}}-{{\alpha }_{1}}-{{\alpha }_{2}}>0$, and ${{\beta }_{2}}+{{\beta }_{3}}-{{\beta }_{1}}>0$ ($E_{13}(\cdot)$ function was studied in \cite{KTU25}).

{\bf Proof. }
The proof of this theorem is similar to the proof of Theorem 2.

According to the integral representation of $E_2$ (see Lemma 2 in \cite{KKT26}), if $\alpha=1$, $0<\beta<1$, $a<0$, $\delta<0$ and $\gamma=\beta$, then the following inequality holds for the bivariate Mittag-Leffler type function: 
\begin{equation}\label{E2>0}
{{E}_{2}}\left( \left. \begin{aligned}
  & \gamma ,1,\gamma ;0,1; \\ 
 & \beta ,\alpha ,\beta +1;\gamma ,\gamma ;1,1; \\ 
\end{aligned} \right|\begin{matrix}
   a{{t}^{\beta }}  \\
   \delta {{t}^{\alpha }}  \\
\end{matrix} \right)>0.    
\end{equation}
From this statement,
$$B={{M}_{1}}\left( x,x \right)=\Gamma \left( \gamma  \right)\int\limits_{0}^{q}{M}\left( t \right){{t}^{\beta }}{{E}_{2}}\left( \left. \begin{aligned}
  & \gamma ,1,\gamma ;0,1; \\ 
 & \beta ,\alpha ,\beta ;\gamma ,\gamma ;1,1; \\ 
\end{aligned} \right|\begin{matrix}
   a{{t}^{\beta }}  \\
   \delta {{t}^{\alpha }}  \\
\end{matrix} \right)dt\ne 0.$$

\section*{Author Contributions}

All authors contributed equally to this work.

\section*{Conflict of Interest}

The authors declare no conflict of interest.

\par

\end{document}